\documentclass[titlepage,twoside,12pt]{article}
\usepackage{amssymb}
\usepackage{amsfonts}
\begin{document}

{\bf \Large Magmatic "Quantum-Like" Systems} \\ \\

{\bf Elem\'{e}r E Rosinger} \\ \\
Department of Mathematics \\
and Applied Mathematics \\
University of Pretoria \\
Pretoria \\
0002 South Africa \\
eerosinger@hotmail.com \\ \\

{\bf Abstract} \\

Quantum computation has suggested, among others, the consideration of "non-quantum" systems which in certain respects may behave "quantum-like". Here, what algebraically appears to be the most general possible known setup, namely, of {\it magmas} is used in order to construct "quantum-like" systems. The resulting magmatic composition of systems has as a well known particular case the tensor products. \\ \\

{\bf 1. Preliminaries} \\

The way usual "non-quantum" systems are composed, among them registers in usual electronic digital computers is the Cartesian product. The "quantum-like" composition, as in the registers of quantum computers, for instance, is done through tensor products. Here we introduce {\it magmatic products} which - algebraically - are the most general extensions of tensor products. Recently, alternative generalizations of tensor products, mostly beyond the confines of algebra, were presented in [3-7]. \\  

We briefly recall the algebraic concept of magma and some of its basic features used in the sequel, [1]. A {\it magma} is
an algebraic structure $( E, \alpha )$ on an arbitrary nonvoid set $E$, where \\

(1.1)~~~ $ \alpha : E \times E \longrightarrow E $ \\

is any binary operation on $E$. \\

Thus in particular, the binary operation $\alpha$ is {\it not} required to be associative, commutative, or to have any
other specific property. In this way, all the usual binary operations in algebra, among them addition in semigroups,
groups, vector spaces and algebras, or multiplication in algebras are in fact magmas. \\

One of the striking features of magmas, which make them so general, is that they are not required to be associative. This however, is a most basic and
common feature of the way usual human languages are constructed and used, even if we are seldom consciously enough aware of it. Indeed, when it comes to
the way meaning is associated to a succession of words in such languages, word successions are on occasion, and obviously, not associative, as simple
examples illustrate it. Often, appropriate punctuation is employed in order to highlight such a lack of associativity, and help specifying the correct
meaning of the respective succession of words. Also, frequently the presence of non-associativity is simply left unmarked, and thus to the ability of the
reader or listener to properly interpret the text. \\

Given, therefore, this ultimate generality, yet at the same time, manifest ubiquity of magmas, it may be natural to enquire to what extent they could be
useful in the study of "quantum-like" behaviour. \\

In this regard, one can note that, according to present understanding, a quantum computer has as an essential component a so called "quantum register"
which is a collection of a number of independent qubits, a number of at least several hundred, in order to have sufficient practical relevance. It
follows that of crucial importance is

\begin{itemize}

\item to set up an effective physical "quantum register" as a {\it composite} quantum system, that is, a composite system with "quantum-like" behaviour,

\item to understand the realms of mathematics underlying the functioning of such a composite system.

\end{itemize}

Indeed, two of the first and basic features which distinguish quantum computers from the usual electronic digital ones are that

\begin{itemize}

\item instead of bits, they operate with qubits which contain far more information,

\item the qubits are composed into a "quantum register" not according to the classical way, that is, through a Cartesian product, but in "quantum-like"
manner, namely, through tensor product.

\end{itemize}

Here we can note a recent fortunate trend in quantum foundational studies, a trend which is engaged in studying larger and deeper mathematical structures than those which happen to be involved actually in the original von Neumann model, see for instance [2] and the literature cited there. And one specific version of this trend is the study of a possible variety of mathematical alternatives to the customary "quantum-like" manner of composing "quantum-like" systems. \\
This, therefore, can be seen as one motivation of the present paper. \\

Returning now to magmas, let us recall the basic concepts and result needed in the sequel. \\

First, about notation which is important in view of the fact that in a magma $( E, \alpha )$ in (1.1) the binary operation
$\alpha$ need not be commutative or associative, therefore, algebraic expressions involving it may lead to a cumbersome notation. \\

Here we introduce certain useful simplifications through the appropriate use of brackets "$($" and "$)$", as well as the
elimination of comas "$,$" and even of the symbol $\alpha$ of the binary operation defining the magma $( E, \alpha )$ in (1.1).
This will deal with the typical non-associativity of the binary operation $\alpha$. \\
Dealing with possible non-commutativity is simple and well known, namely, one does not permute anything at all in a string
of symbols. \\

Let us now indicated the details of the mentioned simplification. \\
Let $a, b \in E$. Then we shall simply denote $\alpha ( a, b ) = a b$, where the order of the arguments $a$ and $b$ is important, since $\alpha$ need
not be commutative. \\
When it comes to three arguments $a, b, c \in E$, then the judicious use of brackets becomes necessary. Indeed, here
enters the fact that $\alpha$ need not be associative, thus we can have $\alpha ( a, \alpha ( b, c ) )$ or $\alpha
( \alpha ( a, b ), c )$, and in general, these two expressions need not be equal. \\

It is clear however that, without any ambiguity, one can for simplicity denote $\alpha ( a, \alpha ( b, c ) ) =
( a ( b c ) )$, and $\alpha ( \alpha ( a, b ), c ) = ( ( a b ) c )$. In other words, one can simply omit $\alpha$, as well
as the comas. \\

For four or more arguments the extension of this simplified notation is obvious. \\

Before going further, we have to note that, as seen below, Cartesian products will be involved in the definition of the
{\it free magma} $M_X$ generated by an arbitrary nonvoid set $X$. \\

Further, we also have to recall that the operation of Cartesian product, when considered rigorously in its notation, is in
fact {\it neither} commutative, {\it nor} associative. \\
Indeed, if $E$ and $F$ are two arbitrary nonvoid sets, then the elements of their Cartesian product $E \times F$ are
denoted by $( x, y )$, while those of $F \times E$ are $( y, x )$, where $x \in E$ and $y \in F$. \\
Given now three arbitrary nonvoid sets $E, F$ and $G$, the Cartesian product $( E \times F ) \times G$ has therefore
elements of the form $( ( x, y ), z )$, while the elements of $E \times ( F \times G )$ are of the form
$( x, ( y, z ) )$, where $x \in E,~ y \in F,~ z \in G$. \\
What is usually done here is simply to make the identification $( ( x, y ), z ) = ( x, ( y, z ) ) = ( x, y, z )$, by
eliminating suitable brackets, and thus obtaining the associativity of Cartesian products by assuming that $( E \times F )
\times G = E \times ( F \times G ) = E \times F \times G$. \\

On the other hand, since magmas are in general {\it not} supposed to be commutative or associative, we shall haven to keep the original
non-commutativity and non-associativity as is in fact already present in the rigorous notation of Cartesian products. And yet, we still can introduce a
simplification in notation, without imposing by that either commutativity or associativity. Namely, we shall keep the brackets, and only have the comas
eliminated. \\

Clearly, the only condition required in order to make the above elimination of brackets and comas free of ambiguity is that the sets
involved in the Cartesian products be such that any concatenation $w_1 w_2$ of elements $w_1, w_2$ in any of the factor sets be uniquely decomposable
into the constituents $w_1$ and $w_2$. And this condition, which can obviously be satisfied without loss of generality, will be assumed in the
sequel. \\

Given now any nonvoid set $X$, we define \\

(1.2)~~~ $ M_X = \bigcup_{1 \leq n < \infty} X_n $ \\

where \\

(1.3)~~~ $ \begin{array}{l}
                    X_1 = X \\
                    \vdots \\
                    X_n = \bigcup_{p+q=n} X_p \times X_q
            \end{array} $ \\ \\

and define on $M_X$ the binary operation \\

(1.4)~~~ $ \gamma : M_X \times M_X \ni ( u, v ) \longmapsto \gamma ( u, v ) = ( u v ) \in M_X $ \\

that is, the concatenation of the elements in $M_X$, plus the enclosure of that concatenation in a pair of brackets. \\

For instance, as defined in (1.2), the first constituent of $M_X$, according to (1.3), is $X_1 = X$, thus there is no
difficulty in understanding it since its elements are simply $x$, with $x \in X$. The second constituent is $X_2 = X_1
\times X_1 = X \times X$, thus with the above simplifying convention about the notation of elements in Cartesian products,
$X_2$ is the set of all elements $( x y )$, with $x, y \in X$, where commutativity is not assumed, hence no particular difficulty either. \\
The specific magmatic feature of $M_X$ starts however with $X_3$, since it relates to the {\it lack} of associativity of
$M_X$, and thus by necessity, it has to refer to at least three elements in $X$. Namely, according to (1.3), we have \\

(1.5)~~~ $ X_3 = ( X_1 \times X_2 ) \cup ( X_2 \times X_1 ) $ \\

which means that $X_3$ is the set of all elements of the form $( x ( y z ) )$ in $X_1 \times X_2$, or $( ( x y ) z )$ in
$X_2 \times X_1$, where $x, y, z \in X$. And clearly, it is {\it not} assumed that $( x ( y z ) ) = ( ( x y ) z )$, {\it
not} even in the case when $x = y = z$. \\

Similarly, $X_4$ is the set of all elements of one of the following five forms $( x_1 ( x_2 ( x_3 x_4 ) ) ),~
( x_1 ( ( x_2 x_3 ) x_4 ) ),~ ( ( x_1 x_2 ) ( x_3 x_4 ) ),~ ( ( x_1 ( x_2 x_3 ) ) x_4 )$ or $( ( ( x_1 x_2 ) x_3 ) x_4 )$,
where $x_1, x_2, x_3, x_4 \in X$. And so on. \\

Here one can note that, given a magma $( E,\alpha )$, the maximum number of different elements of $E$ one can obtain due to non-associativity and
non-commutativity when applying $\alpha$ a number $n$ of times, that is, to $n + 1$ different arguments $a_1, a_2, \ldots , a_n, a_{n+1} \in E$, is
given by the Catalan number $C_n$, thus it grows quite fast with $n$. \\

{\bf Definition 1.1.} \\

$( M_X, \gamma )$ is called the {\it free magma} generated by the nonvoid set $X$.

\hfill $\Box$ \\

Obviously, we have the injective mapping \\

(1.6)~~~ $ \iota_X : X \ni x \longmapsto X = X_1 \subset M_X $ \\

In view of the above, we obtain \\

{\bf Lemma 1.1.} \\

For $n \geq 3$, the constituents $X_p \times X_q$ defining $X_n$ in (1.3) are {\it pairwise disjoint}.

\hfill $\Box$ \\

This makes clear the way $\gamma$ in (1.4) operates, namely \\

(1.7)~~~ $ X_p \times X_q \ni ( u, v ) \longmapsto \gamma ( u, v ) = ( u v ) \in X_n,~~~ p, q \geq 1,~~ p + q = n $ \\

For instance, given $x, y, z, w \in X$, we have \\

$~~~ \begin{array}{l}
             \gamma ( x, y ) = ( x y ),~~~ x, y \in X_1 = X,~~ ( x y ) \in X_1 \times X_1 = X_2 \\ \\
             \gamma ( x, ( y z ) ) = ( x ( y z ) ),~~~ x \in X_1,~~ ( y z ) \in X_2,~~
                                           ( x ( y z ) ) \in X_1 \times X_2 \subset X_3 \\ \\
             \gamma ( ( x y ), ( z w ) ) = ( ( x y ) ( z w ) ),~~~ ( x y ), ( z w ) \in X_2,~~
                                         ( ( x y ) ( z w ) ) \in X_2 \times X_2 \subset X_4
        \end{array} $ \\

and so on. \\

Given now two magmas $( E, \alpha )$ and $( F, \beta )$, then a mapping $f : E \longrightarrow F$ is called a {\it magma morphism}, if and only if \\

(1.8)~~~ $ f ( \alpha ( a, b ) ) = \beta ( f ( a ), f ( b ) ),~~~ a, b \in E $ \\

or equivalently in simpler notation \\

(1.8$^*$)~~~ $ f ( a b ) = f ( a ) f ( b ),~~~ a, b \in E $ \\

The basic result about magmas is in \\

{\bf Theorem 1.1.} \\

Given any mapping $f : X \longrightarrow E$ where $X$ is a nonvoid set and $( E, \alpha )$ is a magma. Then there exists a unique magma morphism
$F : M_X \longrightarrow E$ such that the diagram commutes \\

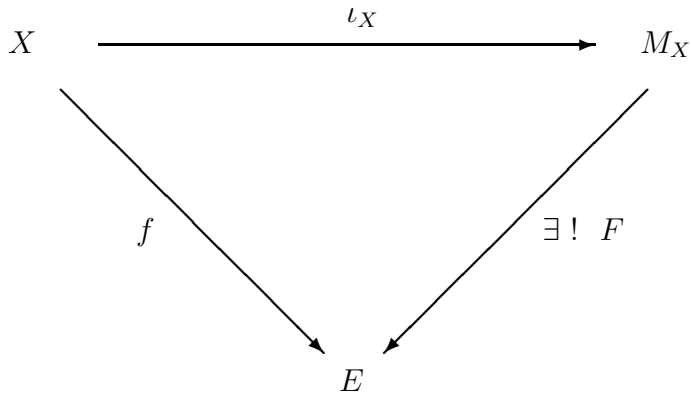
\begin{math}
\setlength{\unitlength}{1cm}
\thicklines
\begin{picture}(13,6)

\put(1,5){$X$}
\put(2.2,5.1){\vector(1,0){6.6}}
\put(9.4,5){$M_X$}
\put(5.5,5.4){$\iota_X$}
\put(1.7,4.5){\vector(1,-1){3.5}}
\put(9.5,4.5){\vector(-1,-1){3.5}}
\put(5.4,0.5){$E$}
\put(2.7,2.5){$f$}
\put(8.1,2.5){$\exists~!~~ F$}

\end{picture}
\end{math}

{\bf Proof.} \\

In view of (1.3), we define $F : M_X \longrightarrow E$ by induction as follows \\

(1.9)~~~ $ X_1 \ni x \longmapsto F ( x ) = f ( x ) \in E $ \\

while for $p, q \geq 1,~~ p + q = n$, we define \\

(1.10)~~~ $ X_p \times X_q \ni ( u, v ) \longmapsto F ( u, v ) = ( F ( u ) F ( v ) ) \in E $ \\

which is a correct definition, due to Lemma 1.1. \\ \\

{\bf 2. "Quantum-Like" Structures} \\

A "quantum-like" register ${\cal R}$, that is a register which is composed "quantum-like" from components ${\cal C}_i$
which have state spaces given by nonvoid
sets $X_i$, with $i \in I$, and are composed in a manner which widely generalizes the usual algebraic composition of
registers in quantum computers, can be
obtained as follows. \\

Let \\

(2.1)~~~ $ X = \prod_{i \in I} X_i $ \\

be the usual Cartesian product of the respective state spaces, thus the state space which would correspond to the usual "non-quantum-like" register with
components ${\cal C}_i$, where $i \in I$. It follows that the elements of $X$ are of the form \\

(2.2)~~~ $ x = ( x_i )_{i \in I} $ \\

where $x_i \in X_i$, with $i \in I$. \\

Then as a {\it first step}, we form the {\it free magma} $( M_X, \gamma )$ generated by $X$. \\

Here, before proceeding with the second step, it is useful further to clarify the notation for the elements of a free magma $M_X$ in general, as defined
in (1.3) for arbitrary nonvoid sets $X$. \\

Namely, given $\xi \in X_n$, with $n \geq 2$, then $\xi$ is a string which consists of $n$ elements $x_1, \ldots , x_n \in
X$ in the given order, and $n - 1$ pairs of matched brackets $(~~~)$. \\

Let us index these bracket pairs by denoting them in some order with $(_1~~~)_1$, \ldots , $(_{n-1}~~~)_{n-1}$. The way these brackets
appear in the string $\xi$ is arbitrary, except for the following restriction, with the order understood from left to
right : \\

(2.3)~~~ If $(_h$ appears before $(_k$, then $)_k$ must appear before $)_h$. \\

In other words, if a pair of brackets $(_h~~~ )_h$ overlaps with a pair of brackets $(_k~~~ )_k$, then one of the pairs must contain the other. \\

Further, given $1 \leq m \leq n$, we denote by $\xi ( m )$ the element $x_m \in X$ in the string $\xi$, as counted from
the left to right. In this way, by discarding now all the brackets, we associate with the string $\xi$, the following
simplified string \\

(2.4)~~~ $ \xi_X = \xi ( 1 ) \ldots \xi ( n ) = x_1 \ldots x_n $ \\

Obviously, the mapping \\

(2.5)~~~ $ \xi \longmapsto \xi_X $ \\

is {\it not} injective when $n \geq 3$. Let us therefore denote \\

(2.6)~~~ $ X_n ( x_1, \ldots , x_n ) = \{~ \xi \in X_n ~~|~~ \xi_X = x_1 \ldots x_n ~\} $ \\

for $x_1, \ldots , x_n \in X$. As follows from an earlier remark, $X_n ( x_1, \ldots , x_n )$ has a number of elements
given by the Catalan number $C_{n-1}$. \\

Clearly, the only condition required in order to make the above free from ambiguity is that the sets $X$ be such that any
concatenation $x_1 x_2 \ldots x_n$ of elements $x_1, x_2, \ldots , x_n \in X$ is uniquely decomposable into the constituents $x_1, x_2, \ldots ,
x_n \in X$. And this condition, which can obviously be satisfied without loss of generality, will be assumed in the sequel. \\

Now, a given $\eta \in X_{n+1}$ is called a {\it replacement} of a certain $\xi \in X_n$, with $n \geq 1$, if and only if, for suitable $1 \leq
h \leq n$ and $x\,'_h, x\,''_h \in X$, we have \\

(2.7)~~~ $ \xi_X = x_1 \ldots x_n,~~~ \eta_X =  x_1 \ldots  x_{h-1} \, x\,'_h \, x\,''_h \, x_{h+1} \ldots  x_n $ \\

and in addition \\

(2.8)~~~ All the $n - 1$ bracket pairs in $\xi$ are preserved in $\eta$, which will \\
         \hspace*{1.3cm} have one more bracket pair so that $\eta$  satisfies condition (2.3). \\

In other words, $x_h$ in $\xi$ was replaced with $x\,'_h \, x\,''_h$, thus leading to $\eta$. However, the additional bracket pair in $\eta$ can be
placed in more that one way, as illustrated in the following examples. Let $\xi = ( x y )$ and let us replace $x$ with $x\,' x\,''$. Then we can have
$\eta = ( x\,' ( x\,'' y ) )$, or $\eta = ( ( x\,' x\,'' ) y )$. Or let $\xi = ( x ( y z ) )$, then we can have $\eta = ( x\,' ( x\,'' ( y z ) ) )$ or
$\eta = ( ( x\,' x\,'' ) ( y z ) )$. \\

We note in this regard that a given element $x \in X$ in the string $\xi$ can be in one of the following situations only, with respect
to its immediate neighbours \\

(2.9)~~~ $ \begin{array}{l}
                    1)~~~ ( x ( \\
                    2)~~~ ) x ) \\
                    3)~~~ ( x y ) \\
                    4)~~~ ( y x )
             \end{array} $ \\

Consequently, when replacing $x$ with $x\,' x\,''$, we can only have the following situations \\

(2.10)~~~ $ \begin{array}{l}
                    1)~~~ ( x\,' (_*\, x\,'' ( \\
                    2)~~~ ( (_*\, x\,' x\,'' )_*\, ( \\
                    3)~~~ ) x\,' )_*\, x\,'' ) \\
                    4)~~~ ) (_*\, x\,' x\,'' )_*\, ) \\
                    5)~~~ ( x\,' (_*\, x\,'' y )_*\, ) \\
                    6)~~~ ( (_*\, x\,' x\,'' )_*\, y ) \\
                    7)~~~ ( y (_*\, x\,' x\,'' )_*\, ) \\
                    8)~~~ ( (_*\, y x\,' )_*\, x\,'' )
            \end{array} $ \\

where the brackets $(_*$ and $)_*$ are all, or part of the new brackets added. \\

Therefore, we denote by \\

(2.11)~~~ $ X_{n+1,\, \xi}\, ( x_1, \ldots , x_{h-1}, x\,'_h, x\,''_h, x_{h+1}, \ldots , x_n ) $ \\

the set of all $\eta \in X_{n+1} ( x_1, \ldots , x_{h-1}, x\,'_h, x\,''_h, x_{h+1}, \ldots , x_n )$ which are a replacement
of $\xi$. \\

Now as a {\it second step}, and assuming arbitrary magma structures $( X_i, \alpha_i )$ on the state spaces $X_i$ of the
components ${\cal C}_i$, with $i \in I$, we define on $M_X$ the equivalence relation $\approx$ as follows. Given $\xi,
\eta \in M_X$, then $\xi \approx \eta$, if and only if $\xi = \eta$, or $\xi$ and $\eta$ can be obtained from one another
by a finite number of the following kind of replacements. \\

Let $\xi \in X_n \subset M_X$, with $n \geq 1$, and assume that \\

(2.12)~~~ $ \xi_X = \xi ( 1 ) \ldots \xi ( n ) = x_1 \ldots x_n $ \\

Let $1 \leq h \leq n$, and assume that \\

(2.13)~~~ $ x_h = ( x_{\,h,\,i} )_{i \in I} $ \\

Further, assume that for some $j \in I$, we have \\

(2.14)~~~ $ x_{\,h,\,j} = \alpha_j ( x\,'_{\,h,\,j}, x\,''_{\,h,\,j} ) $ \\

where $x\,'_{\,h,\,j}, x\,''_{\,h,\,j} \in X_j$. \\

Then $\xi$ is replaced with $\eta \in  X_{n+1,\, \xi}\, ( x_1, \ldots , x_{h-1}, x\,'_h, x\,''_h, x_{h+1}, \ldots , x_n )  \subset M_X$, hence \\

(2.15)~~~ $ \eta_X = x_1 \ldots x_{h-1} \, x\,'_h \, x\,''_h \, x_{h+1} \ldots x_n $ \\

while \\

(2.16)~~~ $ x\,'_h = ( x\,'_{\,h,\,i} )_{i \in I},~~~ x\,''_h = ( x\,''_{\,h,\,i} )_{i \in I} $ \\

with \\

(2.17)~~~ $ x\,'_{\,h,\,i} = \begin{array}{|l}
                                   x_{\,h,\,i} ~~\mbox{if}~ i \neq j \\
                                   x\,'_{\,h,\,j} ~~\mbox{if}~ i = j
                             \end{array}~~~~~~~
                             x\,''_{\,h,\,i} = \begin{array}{|l}
                                   x_{\,h,\,i} ~~\mbox{if}~ i \neq j \\
                                   x\,''_{\,h,\,j} ~~\mbox{if}~ i = j
                             \end{array} $ \\

Based on the above, we can introduce \\

{\bf Definition 2.1.} \\

The quotient set \\

(2.18)~~~ $ \boxtimes_{i \in I} X_i = M_X / \approx $ \\

is called the {\it magmatic product} of the family $( X_i, \alpha_i )$, with $i \in I$.

\hfill $\Box$ \\

Let us recapitulate. \\

The free magma $( M_X, \gamma )$ is defined exclusively in terms of the nonvoid sets $X_i$, with $i \in I$. Then in the definition (2.18) of the magmatic
product $\boxtimes_{i \in I} X_i$, use is made of the equivalence relation $\approx$ whose definition (2.12) - (2.17) involves the magma operations
$\alpha_i$ on $X_i$, with $i \in I$. \\

What is left now is to define the magma structure on the above magmatic product $\boxtimes_{i \in I} X_i$. This will be done in a natural manner by
applying the above quotient operation in (2.18) defined by the equivalence relation $\approx$ not only to $M_X$, but also to the binary operation
$\gamma$ on $M_X$. In this way, we shall obtain a binary operation $\delta$ on $\boxtimes_{i \in I} X_i$ resulted from $\gamma$ through a quotient
operation defined by the equivalence relation $\approx$. \\

For that purpose, we recall a related general result. Let $( E, \alpha )$ be any magma. An equivalence relation $\approx$ on $E$ is called a {\it
congruence}, if and only if for $a, b, c \in E$, we have \\

(2.19)~~~ $ a \approx b ~~\Longrightarrow~~ \alpha ( a, c ) \approx \alpha ( b, c),~~~ \alpha ( c, a ) \approx \alpha ( c, b ) $ \\

Given now a congruence $\approx$ on a magma $( E, \alpha )$, let us define the binary operation $\beta$ on the quotient $E / \approx$ as follows \\

(2.20)~~~ $ \beta ( ( a )_\approx, ( b )_\approx ) = ( \alpha ( a, b ) )_\approx,~~~ a, b \in E $ \\

where for $c \in E$, we denote by $( c )_\approx \in E / \approx$ the $\approx$ equivalence class, or coset of $c$ in $E$. Then $\beta$ is well defined.
Indeed, let $a, a\,', b, b\,' \in E$, with $a \approx a\,'$ and $b \approx b\,'$. Then (2.19) gives $\alpha ( a, b ) \approx \alpha ( a\,', b ) \approx
\alpha ( a\,', b\,' )$. \\

Now we note that the equivalence $\approx$ defined in (2.12) - (2.17) on $( M_X, \gamma )$ is a congruence. Indeed, in view of (2.7) - (2.11), if $\eta
\in X_{n+1}$ is a replacement of $\xi \in X_n$, then for every $\chi \in M_X$ it follows that $( \chi \eta )$ is a replacement of $( \chi \xi )$, while
$( \eta \chi )$ is a replacement of $( \xi \chi )$. \\
Therefore, the desired binary operation $\delta$ on $\boxtimes_{i \in I} X_i$ can be defined by the procedure in (2.21). \\

In this way we arrive to \\

{\bf Definition 2.2.} \\

The magma structure on the magmatic product $\boxtimes_{i \in I} X_i$ is given by \\

(2.21)~~~ $ ( \boxtimes_{i \in I} X_i, \delta ) = ( M_X, \gamma ) / \approx $ \\

{\bf Remark 2.1.} \\

1) The magma $( \boxtimes_{i \in I} X_i, \delta )$ in (2.21) is neither commutative, nor associative, even if all its component magmas $( X_i,
\alpha_i )$, with $i \in I$, are commutative and associative, as seen in the examples in section 3 below. \\

2) In view of (1.6), (2.18), we have the mapping \\

(2.22)~~~ $ X = \prod_{i \in I} X_i \ni x = ( x_i )_{i \in I} \longmapsto \boxtimes_{i \in I} x_i = ( x )_\approx \in \boxtimes_{i \in I} X_i $ \\

which in general need not be injective. \\ \\

{\bf 3. Examples} \\

Let us start by noting that the construction in section 2 above can be performed even in the particular case of one single magma $( X, \alpha )$, that
is, when the index set $I$ has one single element. The result in such a simple case will be called {\it self-magmatic product}, and denoted by \\

(3.1)~~~ $ \boxtimes ( X, \alpha ) $ \\

Obviously, the resulting magma on this self-magmatic product is neither commutative, nor associative, even in case $\alpha$ is both commutative and
associative, or $X$ has only one single element. \\

As a second example, let $( X, + )$ and $( Y, + )$ be two commutative semigroups. Then \\

(3.2)~~~ $ ( X, + ) \boxtimes ( Y, + ) = ( X \boxtimes Y, \delta ) = ( M_{X \times Y}, \gamma ) / \approx $ \\

and it is easy to see that $\delta$ is neither commutative, nor associative. \\

Here we can note that the usual tensor products $\bigotimes$ of semigroups are associative and commutative up to isomorphism, even when the semigroups are not commutative, [3-7]. \\

The magmatic products $\boxtimes$, on the other hand, are never associative, due to the essentially non-associative nature of free magmas, and of the way the equivalence relations $\approx$ used in the construction of magmatic products are defined in (2.12) - (2.17). \\

\end{document}